\documentclass[final,a4paper,oneside]{article}

\usepackage{graphicx}
\usepackage{julia_styles}
\usepackage{amsmath}
\usepackage{amsfonts}
\usepackage{amssymb}
\usepackage{subcaption}
\usepackage{color}
\usepackage{epstopdf}
\usepackage[ruled]{algorithm2e}
\usepackage{pgfplots}
\pgfplotsset{title style={at={(0.5,0.75)}}}
\usepackage[left=3.6cm,right=3.6cm]{geometry}
\title{Problem adapted Hierarchical Model Reduction for the Fokker-Planck equation}

\author{Julia Brunken$^{\dagger}$\thanks{Corresponding author} \and Tobias Leibner\thanks{Applied Mathematics Muenster, University of Muenster, Einsteinstr. 62, 48149 Muenster, Germany ({\tt julia.brunken@uni-muenster.de, tobias.leibner@uni-muenster.de, mario.ohlberger@uni-muenster.de, kathrin.smetana@uni-muenster.de}).} \and Mario Ohlberger$^{ \dagger}$ \and Kathrin Smetana$^{\dagger}$}
\date{}

\begin{document}

\maketitle

\begin{abstract}
In this paper we introduce a new hierarchical model reduction framework for the Fokker-Planck equation.
We reduce the dimension of the equation by a truncated basis expansion in the velocity variable, obtaining a hyperbolic system of equations in space and time. 
Unlike former methods like the Legendre moment models, the new framework generates a suitable problem-dependent basis of the reduced velocity space that mimics the shape of the solution in the velocity variable. 
To that end, we adapt the framework of [M. Ohlberger and K. Smetana. \emph{A dimensional reduction approach based on the application
of reduced basis methods in the framework of hierarchical model reduction}. SIAM J. Sci. Comput., 36(2):A714–A736, 2014] and derive initially a parametrized elliptic partial differential equation (PDE) in the velocity variable. Then, we apply ideas of the Reduced Basis method to develop a greedy algorithm that selects the basis from solutions of the parametrized PDE. Numerical experiments demonstrate the potential of this new method.
\end{abstract}



\pagestyle{myheadings}
\thispagestyle{plain}
\markboth{J. BRUNKEN, T. LEIBNER, M. OHLBERGER, AND K. SMETANA}{MODEL REDUCTION FOR THE FOKKER-PLANCK EQUATION}

\section{Introduction}

Kinetic equations that describe particle flow in phase space (consisting of spatial, temporal, and velocity or momentum variables) arise in various fields. Due to their high dimension, model reduction techniques are crucial for the numerical solution. 

In this paper we introduce a hierarchical model reduction (HMR) approach for kinetic equations, in which the reduced basis is tailored to the specific problem. Similar to the (Legendre) moment closure approach \cite{Levermore}, we project on a rank-$m$ tensor approximation space of the form $V_m := \{ \psi(t,x,v) = \sum_{i=1}^{m} p_i(t,x) \Phi_i(v) \}$, which leads to a coupled system of hyperbolic partial differential equations (PDEs)  for the coefficient functions (moments) $p_i(t,x)$.
However, different from the moment closure approach we do not choose the basis in the velocity direction a priori, e.g.\ as monomials or Legendre polynomials, but rather suggest computing them in a problem adapted a posteriori manner. In detail, we adapt the framework introduced in \cite{Smetana14} and derive first a parametrized (integro-) differential equation in $v$-direction, where the parameters account for the unknown behavior of the density $\psi(t,x,v)$ on the space and time variable. As in \cite{Smetana14} the basis functions are then selected from the solutions of the parametrized differential equation with Reduced Basis (RB) methods (see for instance \cite{RoHuPa08}), exploiting their excellent approximation properties  \cite{DePeWo13}.

In this work, we exemplify the approach for the one-dimensional Fokker-Planck equation. A derivation of the equation can be found in \cite{Hensel,Pomraning}. In \cite{Frank07,Schneider} full-, half-, and mixed-moment models for the Fokker-Planck equation with different closures are developed and numerical experiments are given.

There are several other contributions in which model reduction approaches have been applied to hyperbolic equations and more specifically to the Fokker-Planck equation. The proper generalized decomposition (PGD) method has been introduced and demonstrated for one-dimensional finitely extensible nonlinear elastic (FENE) dumbbells in \cite{Ammar06,Ammar06b}. To this end, it is assumed that the density can be well approximated by a sum of tensor products in which each factor depends on a different variable, leading to a nonlinear variational problem for the successive enhancement of the bases. In \cite{Knezevic_Patera_2010} RB methods are considered for the parametrized Fokker-Planck equation corresponding to a suspension of FENE dumbbells in a prescribed extensional flow. For (parametric) transport dominated problems we refer to the recent contribution of \cite{Dahmen12,Dahmen14} for an appropriate functional analytical setting and for instance to \cite{OhlRav13,TaPeQu15,Welper15} for methods where nonlinear approximations are employed to improve the often deteriorated convergence behavior of RB methods for transport dominated problems. 

We finally mention that the idea of a basis expansion ansatz to reduce the dimension of a problem is not restricted to kinetic equations. In the HMR approach (originally introduced in \cite{Vogelius1} for linear, elliptic PDEs, applied in \cite{Perotto10} and \cite{Smetana14}, and extended to nonlinear PDEs in \cite{SmeOhl14}), the two-dimensional domain of an elliptic PDE is decomposed in a dominant direction in which the problem is solved in full range, and a subordinate transverse direction in which a basis expansion is performed. As it has been demonstrated numerically in \cite{Smetana14} that problem adapted basis functions often provide a better approximation than a priori chosen basis functions like Legendre or trigonometrical polynomials \cite{Vogelius1,Perotto10}, we adapt in this article ideas from \cite{Smetana14} to the Fokker-Planck equation.

The remainder of this article is organized as follows: After defining the problem setting in Section \ref{sect:Problem_Setting}, we transfer the HMR approach to the Fokker-Planck equation and obtain a reduced system in Section \ref{sec:deriv_FP_system}. In Section \ref{sec:deriv_pPDE}, we adapt the procedure from \cite{Smetana14}. We first derive a parametrized PDE in the velocity variable and then develop a basis generation algorithm with ideas from the RB method. Finally, we demonstrate numerical experiments in Section \ref{sect:Numerical_Experiments}, where we compare results from our new framework to the Legendre moment models.

\section{Problem Setting}\label{sect:Problem_Setting}
We consider the one-dimensional Fokker-Planck equation as treated in \cite{Schneider}:
\begin{equation} \label{eq:FP}
\d_t \psi(t,x,v) + v \d_x \psi(t,x,v) + \sigma_{\text{a}}(t,x) \psi(t,x,v) = \frac{T(t,x)}{2} \Delta_v \psi(t,x,v) + Q(t,x,v),
\end{equation}
for $t \in [0, t_1]$, $x \in (a,b)$, $v \in (-1,1)$, and with the Laplace-Beltrami operator on the unit sphere $\Delta_v \psi= \d_v\left ( (1-v^2) \d_v \psi \right ).$
We prescribe the ingoing radiation at the spatial boundary
\begin{equation}
\begin{split} \label{eq:FP_BC}
\psi(t,a,v)=\psi_a(t,v) \;\;  \text{for } v > 0, \quad
\psi(t,b,v)=\psi_b(t,v) \;\;  \text{for } v < 0, \\
\end{split}
\end{equation}
and the initial condition
\begin{equation} \label{eq:FP_IC}
\psi(0,x,v)=\psi_0(x,v).
\end{equation}

\section{Hierarchical Model Reduction} \label{sec:deriv_FP_system}
In this section we adapt the HMR approach \cite{Vogelius1,Perotto10,Smetana14} to the Fokker-Planck equation and thus make the ansatz
\begin{equation} \label{eq:basis_ansatz}
\psi(t,x,v)= \sum_{i=1}^m p_i(t,x) \phi_i(v),
\end{equation}
for a set of $L^2$-orthogonal functions $\phi_i(v),\, i = 1, \dots, m \in \N$. 
This corresponds to the Legendre moment approach, where a Galerkin semi-discretization is done with a polynomial basis, while using other expansion functions. These will be referred to as basis of the reduced velocity space. We will establish a parametrized problem to compute this basis with RB techniques in Section \ref{sec:deriv_pPDE}. The whole problem is then solved in the $(t,x)$-space with fixed $\phi_i(v), \, i=1, \dots, m$. To obtain the reduced system in the $(t,x)$-space, we first assume the functions $\phi_i(v),$  $i=1,\dots, m,$ to be given. 

Inserting \eqref{eq:basis_ansatz} into the Fokker-Planck equation \eqref{eq:FP} leads to:
\begin{equation*}
\begin{split}
\d_t &\left (\sum_{i=1}^m p_i(t,x) \phi_i(v) \right )+ v \d_x\left (\sum_{i=1}^m p_i(t,x) \phi_i(v) \right )  + \sigma_a(t,x)\left (\sum_{i=1}^m p_i(t,x) \phi_i(v) \right )\\
 &= \frac{T(t,x)}{2} \Delta_v\left (\sum_{i=1}^m p_i(t,x) \phi_i(v) \right ) + Q(t,x,v).
\end{split}
\end{equation*}
We test with the reduced velocity basis functions $\phi_k(v),\, k = 1,\dots,m$. Using the short hand notation for the $L^2$-product in velocity space $(f_1,f_2)_v := \int_{-1}^1 f_1(v)f_2(v) \dd v$ this reads for each $k=1,\dots,m$, 
\begin{align} 
\label{eq:FP-RB_semi_Galerkin} &\sum_{i=1}^m (\phi_i(v),\phi_k(v))_v \d_t p_i(t,x) 
+ \sum_{i=1}^m (v \phi_i(v), \phi_k(v) )_v \d_x p_i(t,x) \\
\nonumber &+ \sum_{i=1}^m (\phi_i(v),\phi_k(v))_v \sigma_a(t,x) p_i(t,x) 
+ \sum_{i=1}^m \left( (1-v^2) \d_v \phi_i(v) ,\d_v \phi_k(v)\right )_v  \frac{T(t,x)}{2} p_i(t,x)\\
\nonumber&\quad =  ( Q(t,x,v) ,\phi_k(v))_v .
\end{align}
Setting $\fp(x,t) = (p_1(t,x),\dots, p_m(t,x))^T$, we can write \eqref{eq:FP-RB_semi_Galerkin} as 
\begin{equation*}
\fM \d_t \fp(t,x) + \fD \d_x \fp(t,x) + \sigma(t,x) \fM \fp(t,x) + \frac{T(t,x)}{2} \fS \fp(t,x) = \fq(t,x) 
\vspace{-1em}
\end{equation*}
\begin{alignat*}{2} 
\text{with} \qquad &\fM =(M_{ij})\in \R^{m \times m}, &&M_{ji} := (\phi_i(v),\phi_j(v))_v, \\
&\fD =(D_{ij})\in \R^{m \times m}, &&D_{ji}:= (v \phi_i(v),\phi_j(v))_v, \\
&\fS=(S_{ij}) \in \R^{m \times m}, &&S_{ji}:=\left ( (1-v^2) \d_v \phi_i(v), \d_v \phi_j (v)\right )_v, \\
&\fq(t,x)=(q_k(t,x)) \in \R^m, \qquad &&q_k(t,x) := (Q(t,x,v),\phi_k(v))_v.
\end{alignat*}
As we require the functions $\phi_i(v), \,i=1,\dots,m$ to form a basis of the reduced velocity space, the mass matrix $\fM$ is invertible. We multiply the equation with $\fM^{-1}$ and obtain the linear hyperbolic system
\begin{equation} \label{FP_system}
\d_t \fp + \fM^{-1}\fD \d_x \fp + \left (\sigma \fI_{m \times m} + \frac{T}{2} \fM^{-1}\fS \right ) \fp = \fM^{-1} \fq.
\end{equation}

 

To transfer the boundary and initial conditions \eqref{eq:FP_BC}, \eqref{eq:FP_IC} of the Fokker-Planck equation \eqref{eq:FP} to the hyperbolic system, we project them on the reduced space. After generating a reduced basis, we compute the initial condition vector by
$
\ftp^0(x) \in \R^m$, $\tilde p^0_i(x) := (\psi_0(x,v),\phi_i(v) )_v$, 
$\fp^0(x) := \fM^{-1} \ftp^0(x).
$

To obtain a well-posed system, we may prescribe boundary conditions only for the incoming waves, while outgoing particles are determined by the particle distribution in the domain. In the phase space formulation of the Fokker-Planck equation, this is simply done by prescribing boundary conditions only for the incoming part of the velocity variable. In system \eqref{FP_system} the distinction between incoming and outgoing waves is not so obvious. However, we will use an upwind scheme for the discretization of the system that automatically takes into account the wave directions. We thus do not have to separate the boundary conditions into prescribed and not prescribed parts. 
We define an ``incoming boundary function'' for the whole velocity space $(-1,1)$ by setting the outgoing parts to zero:
\begin{equation*}
\tilde\psi_a(t,v) = 
\begin{cases}
 \psi_a (t,v), &v>0 \\
 0, &v\le 0 ,
 \end{cases}
 \qquad\quad
 \tilde\psi_b(t,v) = 
\begin{cases}
 0,  &v > 0 \\ 
 \psi_b (t,v), &v \le 0.
 \end{cases}
\end{equation*}
Then, the incoming boundary conditions for system \eqref{FP_system} are defined similarly to the initial condition as: $ \ftp^a(t) \in \R^m$, $\tilde p^a_i(t) := (\tilde\psi_a(t,v),\phi_i(v) )_v$, $
\fp^a(t) := \fM^{-1} \ftp^a(t)$. Likewise, we define the right boundary $\fp^b(t)$. 

\section{Generation of problem adapted basis functions}
\subsection{Derivation of a parametrized PDE in the velocity variable}
\label{sec:deriv_pPDE}

To construct basis functions $\phi_i(v)$ dependent on the velocity variable tailored to the considered problem, we adapt the procedure proposed in \cite{Smetana14}. This means that we introduce a parametrized problem for a density $\phi(v)$ solely depending on the velocity, where the parameters include the unknown behavior in the space and time variable. First, we assume (only for the derivation of the parametrized problem) that the solution can be written as 
\begin{equation*}
\psi(t,x,v)\approx P(t,x)\phi(v),
\end{equation*}
where $P(t,x)$ reflects the at this point unknown behavior of the full solution in the space and time variable. 
To find an equation for $\phi(v)$ dependent on $P(t,x)$, we test the Fokker-Planck equation \eqref{eq:FP} with $\varphi(t,x,v)=P(t,x)\tvphi(v)$. This reads
\begin{equation} \label{eq:test_velocity_eqn}
\begin{split}
(\d_t &P(t,x),P(t,x))_{t,x}(\phi(v),\tvphi(v))_v + (\d_x P(t,x),P(t,x))_{t,x} (v \phi(v),\tvphi(v))_v \\
 &+ (\sigma_a(t,x) P(t,x),P(t,x))_{t,x} (\phi(v),\tvphi(v))_v\\
&+\left(\frac{T(t,x)}{2}P(t,x),P(t,x) \right )_{t,x}\left ( (1-v^2)\d_v\phi(v),\d_v\tvphi(v) \right )_v \\
&=\big( (Q(t,x,v),P(t,x) )_{t,x},\tvphi(v) \big )_{v}, 
\end{split}
\end{equation}
where we abbreviate $\int_0^{t_1}\int_a^b f_1(t,x)f_2(t,x)\dd x \dd t =: (f_1,f_2)_{t,x}$.

Since the function $P(t,x)$ is not known at this point, we cannot compute the integrals of the scalar products. Instead, we follow \cite{Smetana14} and introduce a quadrature formula defining the evaluations of the functions in the variables $x$ and $t$  as a parameter. 
In detail, we approximate the inner products in the $(t,x)$-space with the quadrature formula
\begin{equation} \label{eq:quadrature}
\begin{split}
\left (u_1,  u_2 \right )_{t,x} = \int_0^{t_1} \int_a^b u_1(t,x) u_2(t,x) \dd x \dd t  
\approx \sum_{q=1}^{\tilde q} \omega_q u_1(t^q,x^q) u_2(t^q,x^q)
=: \left (u_1 ,  u_2 \right )^q_{t,x},
\end{split}
\end{equation}
where $\omega_q, \, q=1,\dots,\tilde q,$ are the weights and $(t^q,x^q),\, q=1,\dots, \tilde q,$ are the quadrature points. 
The source term is also approximated with the quadrature formula, and we obtain
\begin{equation*}
\begin{split}
((Q,P)_{t,x} , \tvphi)_{v} &= \int_{-1}^{1} \int_0^{t_1} \int_a^b  Q(t,x,v) P(x,t) \tvphi(v)\dd x \dd t  \dd v \\
&\approx \int_{-1}^{1} \left ( \sum_{q=1}^{\tilde q} \omega_q P(t^q,x^q) Q(t^q,x^q,v)\right ) \tvphi(v) \dd v.
\end{split}
\end{equation*}
We replace $(.,.)_{t,x}$ by $(.,.)_{t,x}^q$ and the source term by the source term with quadrature in \eqref{eq:test_velocity_eqn}. 
Then we have an elliptic PDE on the domain $(-1,1)$ that depends on the quadrature points and the evaluations of $P, \d_x P$, and $\d_t P$ at the quadrature points.  
We thus define the parameter vector $\mu$ containing all these values. Since we do not have any information about the behavior of $\phi(v)$ at the boundaries, we impose inhomogeneous Dirichlet boundary conditions and add the two boundary values $\phi(-1)=\phi_l$ and $\phi(1)=\phi_r$ to the parameter. 
The parameter thus has the form 
\begin{equation*}
\mu=\left( (t^q,x^q)_{1 \le q \le \tilde q} , \big (P(t^q,x^q) \big )_{1 \le q \le \tilde q},  \big(\d_x P(t^q,x^q)\big )_{1 \le q \le \tilde q},\big ( \d_t P(t^q,x^q)\big )_{1 \le q \le \tilde q}, \phi_l, \phi_r  \right ).
\end{equation*}
The parameter space, which contains all admissible parameter values of $\mu$, is defined as $\mathcal{D}:=[0,t_1]^{\tilde q} \times [a,b]^{\tilde q} \times I_0^{\tilde q} \times I_1^{\tilde q} \times I_2^{\tilde q} \times I_B^2 \subset \R^P$. The intervals $I_0, I_1, I_2,$ and $I_B$ contain the ranges of $P, \d_x P,\d_t P,$ and the boundary values. They can be chosen using a priori information on the solution or adaptively during the basis construction process using information from the reduced approximations.
Defining a lifting function $\phi_b(v)$ 
 with $\phi_b(-1)=\phi_l$ and $\phi_b(1)=\phi_r$, the parametrized one-dimensional PDE in $v$-direction is: 

Given any $\mu \in \mathcal{D}$, find $\phi(v) \in  H^1((-1,1))$ such that $\phi(v)-\phi_b  \in H^1_0((-1,1))$ and
\begin{align} \label{eq:1DPPDE}
\begin{split}
&a(\mu)\left ( (1-v^2)\d_v\phi(v),\d_v\tvphi(v) \right )_v 
+ b(\mu) (v \phi(v),\tvphi(v))_v 
+ c(\mu) (\phi(v),\tvphi(v))_v\\
&\hspace{5.5cm}= (\hat Q(v,\mu),\tvphi(v))_v \quad\forall \tvphi \in H^1_0((-1,1)),
\end{split} \\
\begin{split}
&\text{where }a(\mu):=\left(\frac{T}{2}P,P \right )^q_{t,x}, \quad
b(\mu):=(\d_x P,P)^q_{t,x},\\ 
&\hspace{1cm}c(\mu):= (\d_t P,P)^q_{t,x} +   (\sigma_a P,P)^q_{t,x},\quad
\hat Q(v,\mu) := \sum_{q=1}^{\tilde q} \omega_q P(t^q,x^q) Q(t^q,x^q,v).
\end{split} \nonumber
\end{align}

\subsection{Basis generation with RB-techniques} \label{sect:Basis_generation}

\begin{algorithm}[t]
 \KwData{$\Xi=\{\phi^n : n=1, \dots, n_{\text{sample}}\},h_{x},m_{\text{max}}$}
 Initialize $\Phi_0=\emptyset$\;
 \For{$m=1:m_{\text{\emph{max}}}$}{
    
      \For{$n=1 :  n_{\text{\emph{sample}}}$}{
                 	     	      
												$\tilde\Phi^n :=\text{GramSchmidtProcess}			(\Phi_{m-1}\cup \{\phi^n\})$\;
						\If{$\dim (\tilde\Phi^n) = m$			}{			
         $\psi^n:=\text{SolveFPSystem} (\tilde\Phi,h_{x})$\;
         	     	      
          $e_n:=\text{L1Error}(\psi^n)$\;
          }
       }
       $\tilde n := \text{argmin}_{n}(e_n)$\; 
    $\Phi_m := \tilde\Phi^{\tilde n}$\; 
 } 
  \caption{GreedyBasisGeneration \label{algorithm_greedy}}
\end{algorithm}

The RB method is a model reduction technique that was derived for the solution of parametrized PDEs. A detailed introduction can be found for instance in \cite{RoHuPa08}.
The method is developed for cases where the solution of the PDE has to be evaluated for many different parameter values, e.g.\ in optimization problems, or where the solution has to be computed very fast for a previously unknown parameter value. 
The main idea is to compute for a judiciously chosen set of parameter values the corresponding solutions to the PDE, the so-called snapshots. Then, the reduced space is chosen as a subspace of the span of these snapshots.
In the standard RB-context, the reduced space is supposed to have good approximation properties for all parameter values $\mu$. To that end, a proper orthogonal decomposition (POD) can be applied to a sufficiently large number of snapshots. Another possibility is a greedy search, where the aim is to iteratively find the parameter value for which the corresponding reduced solution is worst approximated by the basis constucted so far. The snapshot of this parameter is then used to enrich the basis. 

In the HMR context, however, we are interested in a reduced basis that yields a good approximation $\psi(t,x,v)= \sum_{i=1}^m p_i(t,x) \phi_i(v)$ to the full-dimensional reference solution $\psi^{FD}(t,x,v)$. To that end, the snapshots are supposed to resemble the velocity dependence of the full solution for fixed points $(t,x)$, the parameter values only reflect the unknown behavior of the solution. We therefore use a different greedy algorithm than the standard RB-greedy: Instead of iteratively choosing the worst approximated snapshot for the basis enrichment, 
we want to include only the snapshots in the reduced basis that indeed resemble the dependence in velocity direction of the full solution for certain points $(x,t)$. We thus iteratively extend an empty initial basis by choosing for each $m \le m_{\text{max}}$ one ``best approximating'' snapshot to enrich the basis. 
To that end, we define a discretization of system \eqref{FP_system} with mesh size $h_x$ (for details see Section \ref{sect:Implementation}). For each snapshot we compute the solution of system \eqref{FP_system} based on the initial basis extended by the snapshot and estimate the errors. 
Since we did not yet develop an error estimator for the reduced system, we use the real error to the known full solution to find the best approximating snapshot. The whole greedy algorithm that provides a reduced basis of dimension $m_{\text{max}}$ from a snapshot set $\Xi$ that contains of a number of $n_{\text{sample}}$ snapshots is shown in Algorithm \ref{algorithm_greedy}. 

Another possibility for basis generation is a POD. However, an algorithm that simply applies a POD to snapshots from randomly chosen parameter values does not deliver good results. We suppose that this is mainly due to the fact that we do not expect all admissible parameter values to be useful to the full solution, such that snapshots that are not useful to the solution have a too big influence to the resulting basis.
Algorithms that combine greedy and POD strategies by an adaptive choice of suitable snapshots are subject of future work.

\section{Numerical Experiments} \label{sect:Numerical_Experiments}
\subsection{Implementation} \label{sect:Implementation}

Our implementation of the whole framework is realized in \mbox{pyMOR}, a software library in the Python language for model reduction algorithms, especially for the realization of RB methods, see \cite{pyMOR}. 
For the computation of snapshots of the parametrized PDE \eqref{eq:1DPPDE}, the basis generation, and the computation of the system matrices of system \eqref{FP_system} a finite element discretization with mesh size $h_v$ is used. 
For the computation of system \eqref{FP_system} we use forward Euler time stepping, where the time step is adaptively chosen small enough to guarantee a stable system dependent on the basis, together with a finite volume scheme with mesh size $h_x$ and an upwind flux for linear systems.
The source terms are included in the finite volume evaluation by an unsplit method.

After computing the solution of the fully discrete version of system \eqref{FP_system}, we calculate the spatial density  $\psi^{(0)}(t,x) := \int_{-1}^1 \psi(t,x,v) \dd v$ of the solution.
In all numerical tests we compare two solutions by computing the relative $L^1$-error between the respective spatial densities.

We compare the different reduced solutions with a full-dimensional finite difference discretization implemented in the DUNE framework \cite{dune}. The implementation is orientated on the MATLAB-based framework described in \cite{Schneider}. As a time stepping scheme we use the generalized Runge-Kutta solver for stiff problems GRK4T (see \cite{Hairer}).

\subsection{Model Problem}

As a model problem we use the \emph{SourceBeam} test case from \cite{Schneider}. This test case models a beam entering the domain on the left bound. An additional source is located in the middle of the domain. Absorption and transport coefficient are varied through the domain to show different behavior of the solution.

\begin{figure}
\centering
\begin{subfigure}{0.24\textwidth}
\includegraphics[clip=True, trim=2cm 4cm 2cm 1.5cm,width=\linewidth]{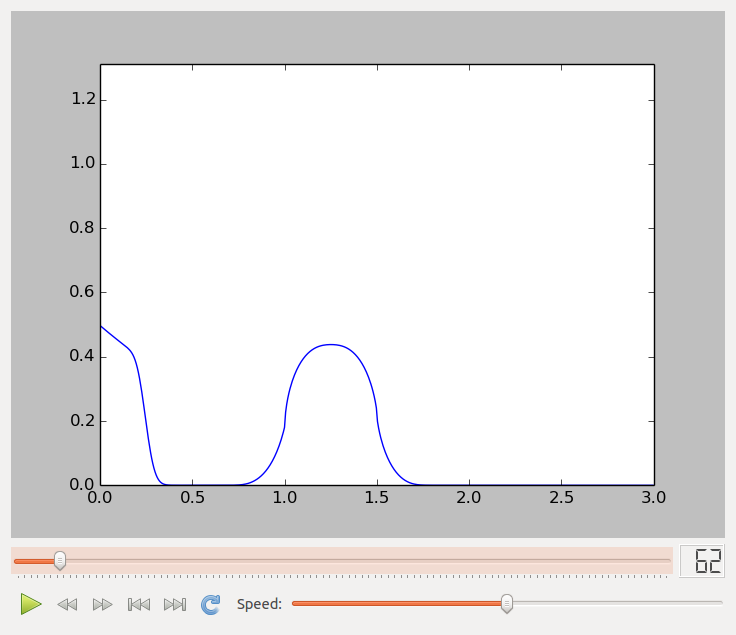}
\caption{$t=0.25$}
\end{subfigure}
\begin{subfigure}{0.24\textwidth}
\includegraphics[clip=True, trim=2cm 4cm 2cm 1.5cm,width=\linewidth]{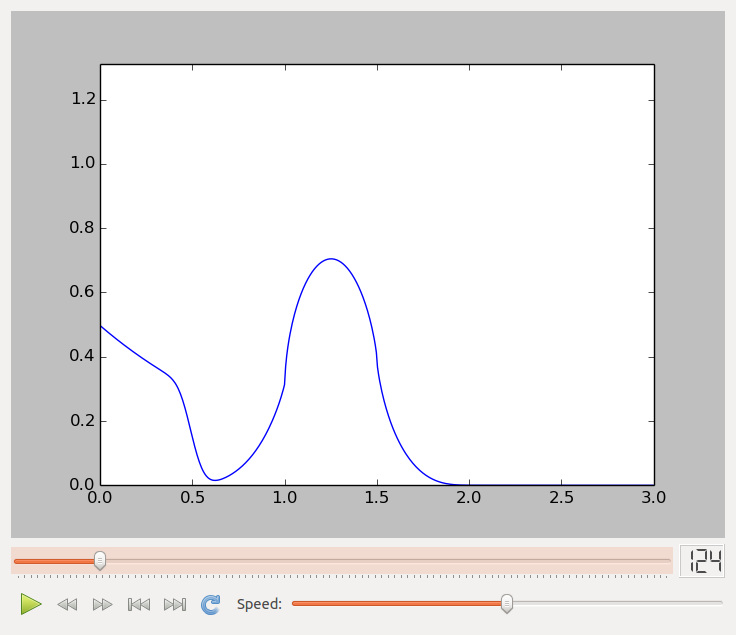}
\caption{$t=0.5$}
\end{subfigure}
\begin{subfigure}{0.24\textwidth}
\includegraphics[clip=True, trim=2cm 4cm 2cm 1.5cm,width=\linewidth]{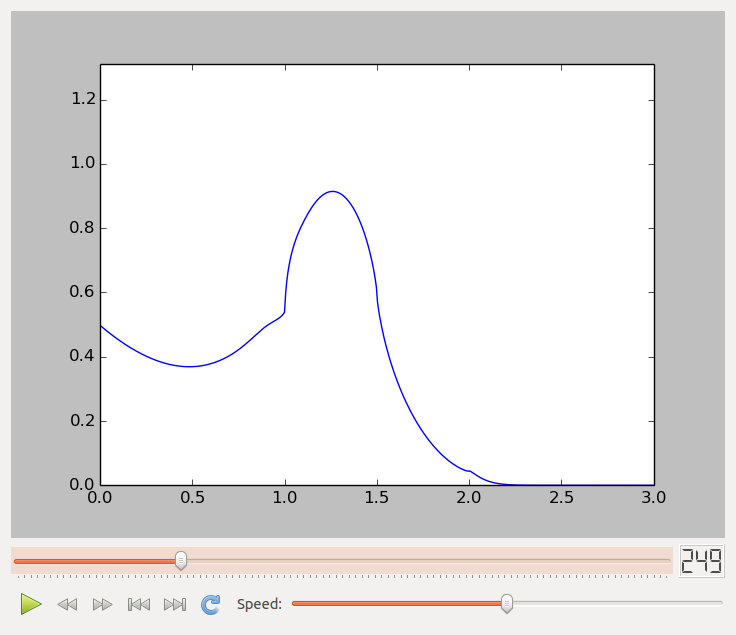}
\caption{$t=1$}
\end{subfigure}
\begin{subfigure}{0.24\textwidth}
\includegraphics[clip=True, trim=2cm 4cm 2cm 1.5cm,width=\linewidth]{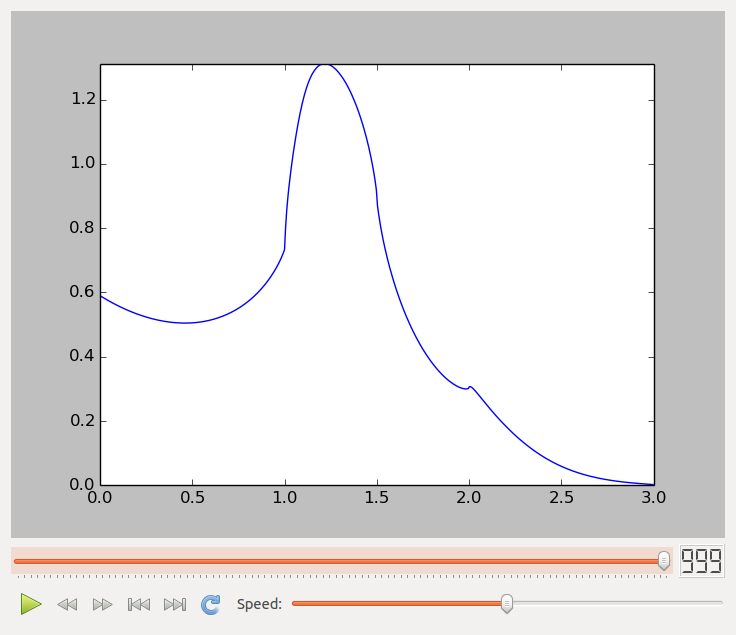}
\caption{$t=4$}
\end{subfigure}
\vspace{-12pt}
\caption{\footnotesize{Density of the full-dimensional finite difference solution for the SourceBeam test case in space for different times.} \label{fig:SourceBeam}}
\vspace{-12pt}
\end{figure}

We consider a domain $x \in (0,3)$ in time $t \in [0,4]$. We define
\begin{equation*} 
\sigma_a(x):= 
\begin{cases}
1, &x \le 2 \\
0, & \text{else,}
\end{cases}
\quad
T(x) :=
\begin{cases}
0, & x \le 1 \\
2, & 1 < x \le 2 \\
10, & \text{else,}
\end{cases}
\quad
Q(x):=
\begin{cases}
1, & 1 \le x \le 1.5 \\
0, & \text{else}.
\end{cases}
\end{equation*}
Initial and boundary conditions are $\psi(0,x,v)=10^{-4}$, $\psi(t,0,v>0)= \delta (v-1)$, $\psi(t,3,v<0)= 10^{-4}$. Plots of the density of the full-dimensional finite difference solution for the SourceBeam test case are shown in Figure \ref{fig:SourceBeam}.

\subsection{Validation of the Basis Generation Process}

As a first validation step and as a benchmark for the whole framework we want to test Algorithm \ref{algorithm_greedy} with ``snapshots'' that indeed \emph{do} resemble the velocity dependence of the full solution for different fixed $(x_i,t_i)$. In this way, we are able to assess the potential of the whole framework when later using snapshots from the parametrized PDE \eqref{eq:1DPPDE} that are \emph{supposed} to resemble the velocity dependence.

To that end, we define ``truth-snapshots'' from a full-dimensional finite difference solution $\psi^{\text{FD}}(t,x,v)$ for points $(t_i,x_i)$ as $s_i(v) := \psi^{\text{FD}}(t_i,x_i,v).$ 
We want to compare different models for different discretization mesh sizes $\tilde h = 2^{-n}, n \in \N$. To that end, for a fixed
mesh size $\tilde h$ we compute the full-dimensional finite difference solution with mesh sizes $h_x=h_v=\tilde h$. Then, the ``truth-snapshots'' are evaluated from this solution at 192 points $(x_i,t_i)$ arranged on a coarse rectangular grid of 12 points in $x$- and 16 points in $t$-direction.
We use Algorithm \ref{algorithm_greedy} with $\Xi=\{s_i(v): i=1, \dots, 192\}$, $m_{\text{max}}=13$ and the mesh size $h_x=\tilde h$ for basis generation.
We compute the solutions of the FP system for different model orders $m \le m_{\text{max}}$ and again with mesh size $h_x=\tilde h$. Then we compare our new model with the Legendre moment models with mesh size $h_x=\tilde h$ implemented in our framework that are derived by simply choosing the first $m$ Legendre polynomials that are discretized with mesh size $h_v=\tilde h$ as reduced velocity basis. 
\begin{figure}
\centering
\footnotesize{
\begin{tikzpicture}
\definecolor{c1}{RGB}{228,26,28}
\definecolor{c2}{RGB}{55,126,184}
\definecolor{c3}{RGB}{77,175,74}
\definecolor{c4}{RGB}{152,78,163}
\definecolor{c5}{RGB}{255,127,0}
\definecolor{c6}{RGB}{166,86,40}
\begin{semilogyaxis}
[width=195pt, ymax=1., ymin=0.005,xmin=1,xmax=13, title=Legendre, xlabel=m,ylabel=$L^1$ error]
\addplot [color=c1,dashed,domain=1:13] {0.0888644214802}; 
\addplot [color=c1,mark=*, mark size=1.5pt] table [col sep=comma]{Legendre_Errors_m3.csv};
\addplot [color=c2,dashed,domain=1:13] {0.0493703043315};
\addplot [color=c2,mark=square, mark size=1.5pt] table [col sep=comma]{Legendre_Errors_m4.csv};
\addplot [color=c3,dashed,domain=1:13] {0.0271248526792};
\addplot [color=c3,mark=triangle*] table [col sep=comma]{Legendre_Errors_m5.csv};
\addplot [color=c4,dashed,domain=1:13] {0.0146163635687};
\addplot [color=c4,mark=diamond] table [col sep=comma]{Legendre_Errors_m6.csv};
\addplot [color=c5,dashed,domain=1:13] {0.00754860362137};
\addplot [color=c5,mark=pentagon*] table [col sep=comma]{Legendre_Errors_m7.csv};
\end{semilogyaxis}
\end{tikzpicture}
\begin{tikzpicture}
\definecolor{c1}{RGB}{228,26,28}
\definecolor{c2}{RGB}{55,126,184}
\definecolor{c3}{RGB}{77,175,74}
\definecolor{c4}{RGB}{152,78,163}
\definecolor{c5}{RGB}{255,127,0}
\begin{semilogyaxis}
[width=195pt, ymax=1, ymin=0.005,xmin=1, xmax=13, legend columns=-1, legend to name=legendname, title=Greedy ``truth snapshots'', xlabel=m]
\addplot [color=c1,dashed,domain=1:13, forget plot] {0.0888644214802}; 
\addplot [color=c1,  mark=*, mark size=1.5pt] table [x index={1}, y index={2}, col sep=comma]{MatlabSnapshotGreedy_m3.csv};
\addlegendentry{$\tilde h=2^{-3}$;}
\addplot [color=c2,dashed,domain=1:13, forget plot] {0.0493703043315};
\addplot [color=c2,mark=square, mark size=1.5pt] table [x index={1}, y index={2}, col sep=comma]{MatlabSnapshotGreedy_m4.csv};
\addlegendentry{$\tilde h=2^{-4}$;}
\addplot [color=c3,dashed,domain=1:13, forget plot] {0.0271248526792};
\addplot [color=c3,mark=triangle*] table [x index={1}, y index={2}, col sep=comma]{MatlabSnapshotGreedy_m5.csv};
\addlegendentry{$\tilde h=2^{-5}$;}
\addplot [color=c4,dashed,domain=1:13, forget plot] {0.0146163635687};
\addplot [color=c4,mark=diamond] table [x index={1}, y index={2}, col sep=comma]{MatlabSnapshotGreedy_m6.csv};
\addlegendentry{$\tilde h=2^{-6}$;}
\addplot [color=c5,dashed,domain=1:13, forget plot] {0.00754860362137};
\addplot [color=c5,mark=pentagon*] table [x index={1}, y index={2}, col sep=comma]{MatlabSnapshotGreedy_m7.csv};
\addlegendentry{$\tilde h=2^{-7}$;}
\end{semilogyaxis}
\end{tikzpicture}
\centering
\begin{tikzpicture}
\definecolor{c1}{RGB}{228,26,28}
\definecolor{c2}{RGB}{55,126,184}
\definecolor{c3}{RGB}{77,175,74}
\definecolor{c4}{RGB}{152,78,163}
\definecolor{c5}{RGB}{255,127,0}
\ref{legendname}
\end{tikzpicture}
\caption{\footnotesize{Errors of reduced models for model orders $m=1,\dots,13$ and mesh sizes $\tilde h = 2^{-3}, \dots, 2^{-7}$ to the reference solution with mesh size $h_x=h_v=2^{-9}$. Left: Legendre moments, right: greedy basis from ``truth-snapshots''. Dashed lines indicate the errors of the full dimensional solutions for the respective mesh sizes $h_x=h_v=2^{-3}, \dots, 2^{-7}$.}}
\label{fig:TruthSnapshots}
\vspace{-10pt}
}
\end{figure}
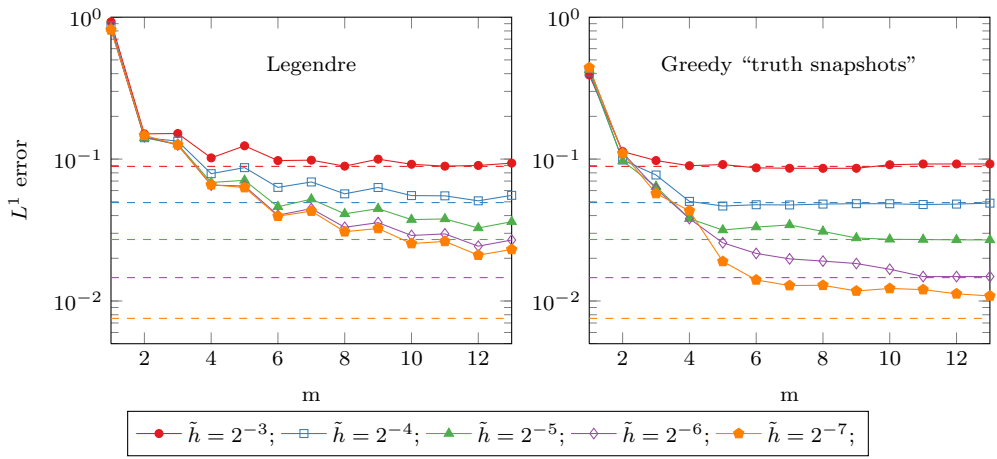

Figure \ref{fig:TruthSnapshots} shows the resulting errors of both methods with mesh sizes $\tilde h = 2^{-3}, \dots, 2^{-7}$ to the reference solution with mesh size $h_x=h_v=2^{-9}$.
To assess the quality of the solutions, the dotted lines indicate the error of the full dimensional finite difference solution with the corresponding mesh sizes $h_x=h_v=2^{-3}, \dots, 2^{-7}$. 

Our new method clearly shows that the Legendre polynomials are not the ideal choice for a basis of the reduced velocity space: 
The greedy basis method provides in all cases a better basis for a given mesh size and model order. The model error rapidly decreases for low model orders and then stagnates about the range of the discretization error that is reached at model orders $m=4$ for $h_x=2^{-3}, 2^{-4}$, $m=9$ for $h_x=2^{-5}$, and $m= 11$ for $h_x=2^{-6}$.
In contrast, the Legendre moment solutions reach the discretization error in the shown model order range $m \le 13$ only for $h_x=2^{-3}, 2^{-4}$. The close proximity of the errors for mesh sizes $h^{-6}$ and $h^{-7}$ indicates a dominance of the model error, whereas the errors of the greedy solutions for $h^{-6}$ and $h^{-7}$ differ for $m \ge 5$, so that we might see a combination of discretization and model error.

These results show the potential of our new method, since the greedy algorithm provided with the ``right'' snapshots shows a better convergence behavior than the Legendre moments.

\subsection{Test with snapshots from the parametrized PDE}

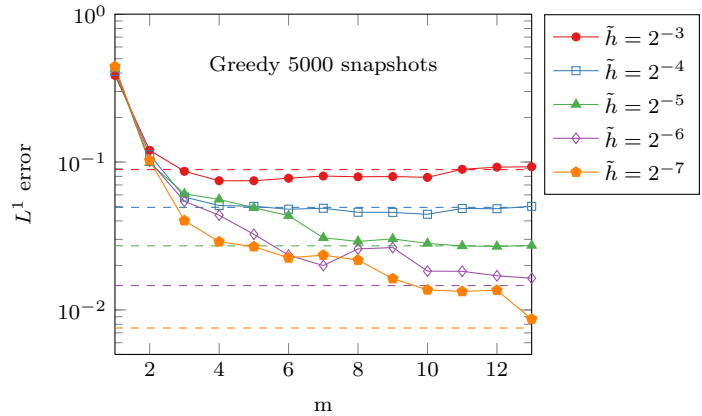
\begin{figure}
\centering
\footnotesize{
	\begin{tikzpicture}
\definecolor{c1}{RGB}{228,26,28}
\definecolor{c2}{RGB}{55,126,184}
\definecolor{c3}{RGB}{77,175,74}
\definecolor{c4}{RGB}{152,78,163}
\definecolor{c5}{RGB}{255,127,0}
\begin{semilogyaxis}
[width=201pt, ymax=1, ymin=0.005,xmin=1, xmax=13,  title=Greedy 5000 snapshots, legend style={legend pos=outer north east}, xlabel=m,ylabel=$L^1$ error]
\addplot [color=c1,dashed,domain=1:13, forget plot] {0.0888644214802}; 
\addplot [color=c1,  mark=*, mark size=1.5pt] table [x index={1}, y index={2}, col sep=comma]{New_Grid_Greedy_snaps=5000,_nx=24.csv};
\addlegendentry{$\tilde h=2^{-3}$}
\addplot [color=c2,dashed,domain=1:13, forget plot] {0.0493703043315};
\addplot [color=c2,mark=square, mark size=1.5pt] table [x index={1}, y index={2}, col sep=comma]{New_Grid_Greedy_snaps=5000,_nx=48.csv};
\addlegendentry{$\tilde h=2^{-4}$}
\addplot [color=c3,dashed,domain=1:13, forget plot] {0.0271248526792};
\addplot [color=c3,mark=triangle*] table [x index={1}, y index={2}, col sep=comma]{New_Grid_Greedy_snaps=5000,_nx=96.csv};
\addlegendentry{$\tilde h=2^{-5}$}
\addplot [color=c4,dashed,domain=1:13, forget plot] {0.0146163635687};
\addplot [color=c4,mark=diamond] table [x index={1}, y index={2}, col sep=comma]{New_Grid_Greedy_snaps=5000,_nx=192.csv};
\addlegendentry{$\tilde h=2^{-6}$}
\addplot [color=c5,dashed,domain=1:13, forget plot] {0.00754860362137};
\addplot [color=c5,mark=pentagon*] table [x index={1}, y index={2}, col sep=comma]{New_Grid_Greedy_snaps=5000,_nx=384.csv};
\addlegendentry{$\tilde h=2^{-7}$}
\end{semilogyaxis}
\end{tikzpicture}
}
\vspace{-6pt}
\caption{\footnotesize{Errors of the greedy method with 5000 snapshots from randomly chosen parameter values to the reference solution for model orders $m=1,\dots,13$ and mesh sizes $\tilde h = 2^{-3}, \dots, 2^{-7}$.  Dashed lines indicate the errors of the full dimensional solutions for the respective mesh sizes $h_x=h_v=2^{-3}, \dots, 2^{-7}$.}}
\label{fig:Greedy5000}
\vspace{-24pt}
\end{figure}

We now want to test our framework with snapshots from the parametrized PDE \eqref{eq:1DPPDE}. We use only one quadrature point $q=(t^1,x^1)$ in the quadrature formula \eqref{eq:quadrature} to limit the dimension of the parameter space. We choose the intervals defining the parameter space as $t^1 \in [0,4]$, $x^1 \in (1,3)$, $P(q) \in [0.01,1.2], \d_x P(q) \in [-5.4,0.9],$ and $\d_t P(q) \in [0,5]$ empirically from already known solutions and the boundary value interval as $[0,1]$. Note that we changed the interval for $x^1$ to assure that $a(\mu)$ is nonzero.
We again want to compare different mesh sizes $\tilde h = h_x = h_v$. To that end, for different mesh sizes $h_v$ we compute $5000$ snapshots from randomly chosen parameter values. We then apply Algorithm \ref{algorithm_greedy} for this snapshot set and for the corresponding mesh size $h_x=h_v$. Figure \ref{fig:Greedy5000} shows the resulting errors to the reference solution. 

For most mesh sizes, the results are similar to the results for the truth-snapshots. For larger mesh sizes, the errors are in parts better than the corresponding discretization errors. This might be due to the different discretization frameworks (full dimensional finite difference and combined finite volume/finite element method).
For the smaller mesh sizes, the truth-snapshot solutions are mostly, but not always, slightly better than the greedy snapshot solutions. However, for all shown mesh sizes the greedy solution errors are at the same range as the discretization error for $m=13$. Additionally, all greedy solutions are better than the Legendre moment models comparing the same model orders and grid sizes. 

We thus conclude that the parametrized PDE provides snapshots that can lead to good approximation properties of the reduced solutions. 
However, the basis generation process in the described form is computationally too expensive to outperform the whole Legendre moment method at the moment.

\section{Conclusions}

In this article we introduced a new hierarchical model reduction approach for the Fokker-Planck equation that consists of a (truncated) basis expansion in the velocity variable. To find a reduced space for this expansion that is directly tailored to the problem, a parametrized PDE in velocity direction is derived and a greedy algorithm for the basis generation inspired by RB methods is developed.
A benchmark reduced model and a numerical experiment show good convergence behavior and thus the potential of our new method. The derivation of an a posteriori error estimator and the improvement of the basis generation algorithm are subject to future work.

\section*{Acknowledgement} 
{\it 
The authors would like to thank Florian Schneider for fruitful discussions and for providing us 
with his code on moment methods for kinetic equations, 
which we used as an orientation for our own implementation of both the hyperbolic system and the reference solution.
}

\bibliographystyle{abbrv}
\bibliography{literature} 

\end{document}